**Four Integer Factorization Algorithms**
**N. A. Carella, September, 2010**

*Abstract:* The theoretical aspects of four integer factorization algorithms are discussed in detail in this note. The focus is on the performances of these algorithms on the subset of difficult to factor balanced integers $N = pq$, $p < q < 2p$. The running time complexity of these algorithms ranges from deterministic exponential time complexity $O(N^{1/2})$ to heuristic and unconditional logarithmic time complexity $O(\log(N)^c)$, $c > 0$ constant.

## 1. INTRODUCTION

Let $p$ and $q$ be a pair of primes. As far as time complexity is concerned, the subset of balanced integers $N = pq$, where $p < q < ap$, and $a > 0$ is a fixed parameter, is the most important subset of integers.

The subset of balanced integers $B(x) = \{ N = pq \le x : p < q < ap \}$ of cardinality $B(x) = O(x/\log^2 x)$ has zero density in the set of all nonnegative integers, see Proposition 4. Accordingly, the factorization of a random integer is unlikely to be as difficult as a balanced integer of the same size.

This article discusses the theoretical aspects of four integer factorization algorithms acting on the subset of balanced integers in details. These algorithms are described in Theorems 2, 3, 5, 13 and 14, respectively. The emphasis is on the performances of these algorithms on the subset of difficult to factor balanced integers. The running time complexity of these algorithms ranges from deterministic exponential time complexity $O(N^{1/2})$ to heuristic and unconditional deterministic logarithmic time complexity $O(\log(N)^c)$, $c > 0$ constant. The innovation here involves a technique for generating systems of polynomials equations for the integer factorization problem by means of elementary polynomial operations and nonelementary polynomial operations, see Theorem 13 and 14.

The standard references [CE], [CP], [LA], [RL], [SP], [WS] and others, provide extensive details on the theory of integer factorizations.

## 2. DETERMINISTIC EXPONENTIAL TIME $O(N^{1/2})$

Any integer $N \ge 1$ has a representation as a difference of two squares

$$4N = x^2 - y^2, \tag{1}$$

and any solution of (1) is of the form $x = p + q$, $y = q - p$, where $p$, and $q$ are factors of $N$, but not necessarily primes. The extreme solution $x = N + 1$, $y = N - 1$ does not lead to a nontrivial factorization of $N$, so it is viewed as the trivial solution. The factors of an arbitrary integer $N$, which can be composites or primes, vary from $p = q = N^{1/2}$, to $p = N/2$, $q = 2$, and $p = N$, $q = 1$. Equality occurs if and only if $N$ is a square.

A prime number $N = x^2 - y^2 > 2$ has a unique solution of large integers $x > N^{1/2}$, $y \geq 1$, often consecutive. This solution is also viewed as a trivial solution since it does not lead to a factorization. But if $N$ is not a prime, then there is a nontrivial solution such that $x = x_k \geq N^{1/2}$ is an integer in the sequence of integers

$$x_0 = \sqrt{N}, \; x_1 = \sqrt{N} + 1, \; x_2 = \sqrt{N} + 2, \; ..., \; x_n = (N+1)/2 \,. \tag{2}$$

Technically $x_k = [\sqrt{N}] + k$ or $x_k = [\sqrt{N} + k]$, where the bracket $[\,x\,]$ defines the largest integer function, however the bracket is often omitted to simplify the notations.

***Proposition* 1.** (i) An integer $N$ has a representation as $N = x^2 - y^2$ if and only if $N \neq 4M \pm 2$.
(ii) The number of solutions is $O(\log N)$ on average, and at most $o(N^\epsilon)$ solutions, $\epsilon > 0$.
(iii) A prime $N > 2$ has a unique representation as a difference of square integers.

Proof: A simple congruence verification shows that $4M \pm 2 = x^2 - y^2$ has no solutions for any $M \geq 0$. Further, it is known that almost every integer has less than $2^{2\log\log(N)}$ divisors, and at most $2^{2\log(N)/\log\log(N)}$ divisors, see [HW, p. 345]. The claims follow from these observations. ∎

The *difference of squares method* effectively handles any integer $N \geq 1$ with a pair of factors in the range

$$\sqrt{N} - c_0 N^{1/4} \log^c N < p < \sqrt{N} < q < \sqrt{N} + c_0 N^{1/4} \log^c N \,,$$

with $c_0$, $c >$ constants. The basic structure of this method is attributed to Fermat, but it is not clear if the time complexity analysis was known before modern time.

***Theorem* 2.** (Fermat) If the factors of an integer $N = pq$ satisfy $|\sqrt{N} - p| = O(N^{1/4} \log^c N)$ and $|\sqrt{N} - q| = O(N^{1/4} \log^c N)$, then it can be decomposed in deterministic logarithmic time complexity $O(\log(N)^c)$, $c > 0$ constant.

Proof: Write the difference of squares equation as $y^2 = x^2 - 4N$, and consider the finite sequence of integers $x_k = 2\sqrt{N} + k \leq p + q$, $k \geq 0$. In addition, suppose that $q - p \leq 2N^\alpha \log^{c/2} N$, $\alpha \geq 0$. Then

$$(2\sqrt{N} + k)^2 - 4N = (q - p)^2 \leq 4N^{2\alpha} \log^c N \tag{3}$$

for some integers $k \geq 0$. Expanding the left side and simplifying it, quickly lead to the inequalities $0 \leq k \leq 2N^{2\alpha - 1/2} \log^{c/2} N$. Therefore, if the parameter $\alpha = 1/4$, then the difference $q - p = O(N^{1/4} \log^{c/2} N)$ and $0 \leq k \leq \log^c N$, so the algorithm runs in deterministic logarithmic time complexity $O(\log(N)^c)$, $c > 0$ constant.

∎

**Note.** The standard term *polynomial time* has been replaced with the more descriptive term *logarithmic time*. This is patterned after the closely related term *exponential time*.

Extending the domain of the algorithm (Theorem 2) to $\sqrt{N} - O(N^{1/2}) < p < \sqrt{N} < q < \sqrt{N} + O(N^{1/2})$ turns this technique into a general purpose integer factorization algorithm. This is accomplished by simply continuing the search in (2) for a nontrivial solution of $4N = x^2 - y^2$. The generalized algorithm has deterministic exponential time complexity $O(N^{1/2})$.

Many techniques have been devised to expedite the running time complexity of this algorithm (Theorem 2). A few of these optimizing techniques are developed in [LN], [MK] and others.

***Example*** 1. There is a small subset of balanced integers $N = pq$, $p < q < 2p$, with $p + q = m(m + 1)/2 + E$ or $p + q = m^2 + E$, where $E = O(N^\varepsilon)$ is a small error. This subset of balanced integers has deterministic exponential time complexity $O(N^{1/4+\varepsilon})$, $\varepsilon > 0$. In particular, the case $p + q = m(m + 1)/2$ can be handled with the finite sequence of integers

$$
\begin{aligned}
x_0^2 &= \left(m(m+1)/2\right)^2, & y_0^2 &= x_0^2 - 4N, \\
x_1^2 &= x_0^2 + (m+1)^3, & y_1^2 &= x_1^2 - 4N, \\
x_2^2 &= x_1^2 + (m+2)^3, & y_2^2 &= x_2^2 - 4N, \\
&\cdots & &\cdots \\
x_k^2 &= x_{k-1}^2 + (m+k)^3, & y_k^2 &= x_k^2 - 4N,
\end{aligned}
\tag{4}
$$

where $m \geq [2N^{1/4}]$, and the $k$th integer $x_k^2 = x_{k-1}^2 + (m+k)^3 = 1 + 2^3 + \cdots + (m+k)^3$. This acceleration technique works whenever the sum $p + q = m(m + 1)/2$ (or the difference $q - p$) is a triangular number.

Assuming the constraints $\sqrt{.5N} < p < \sqrt{N} < q < \sqrt{2N}$, the estimated starting point and the stopping point are deduced from the inequalities $1.5\sqrt{N} \leq p + q = m(m+1)/2 \leq 3\sqrt{N}$. In particular, these inequalities imply that a brute search of at most $2N^{1/4}$ integers is needed to determine a nontrivial solution of $4N = x^2 - y^2$.

For a concrete instance, take $N = 193933249$, (this corresponds to a random triangular sum $p + q$ of 5 digits) and put $m = [2N^{1/4}] = 236$. The appropriate finite sequence of integers { $x_k : k \geq 0$ } is

$$
\begin{aligned}
x_0^2 &= \left(m(m+1)/2\right)^2 = 782097156 = 27966^2, & y_0^2 &= x_0^2 - 4N = 6364160 \neq \text{perfect square}, \\
x_1^2 &= x_0^2 + (m+1)^3 = 795409209 = 28203^2, & y_1^2 &= x_1^2 - 4N = 19676213 \neq \text{perfect square}, \\
x_2^2 &= x_1^2 + (m+2)^3 = 808890481 = 28441^2, & y_2^2 &= x_2^2 - 4N = 33157485 \ \neq \text{perfect square}, \\
&\cdots & &\cdots \\
x_k^2 &= x_{k-1}^2 + (m+k)^3 = 893412100 = 29890^2, & y_k^2 &= x_k^2 - 4N = 117679104 = 10848^2.
\end{aligned}
\tag{5}
$$

Accordingly, $4N = 4 \cdot 193933249 = 29890^2 - 10848^2$. This requires just $k = 9$ steps. In contrast, the standard Fermat difference of squares algorithm requires about $k = 2039$ steps starting at $x_k = [2\sqrt{N}] + k$, $k \geq 0$.

For each fixed triangular number $m(m + 1)/2$ in the range $2\sqrt{N} \leq m(m+1)/2 \leq 3\sqrt{N}$, there are about $2N^{1/2-\varepsilon}$ prime pairs $p$, $q$, such that $p + q = m(m + 1)/2$. Consequently, the subset of balanced integers $N = pq$ with triangular sums is small, its cardinality is approximately $O(N^{3/4-\varepsilon})$, $\varepsilon > 0$. In general, the number of unrestricted representations of an integer $n \in \mathbb{N}$ as a sum $n = p_1 + p_2 + \cdots + p_k$ of $k$ primes has the asymptotic formula $R_k(n) = (cn^{k-1}/(k-1)!\log^k n)(1 + O(1/\log n))$.

It is quite plausible that given $u_1$, which is the most significant $(\log N)/4$ bits of $x = N^{1/2} + u_1 N^{1/4} + u_0 \leq (p+q)/2$, where $0 \leq |u_0|$, $|u_1| < N^{1/4}$, and $u_0$ is unknown, the difference of squares method can factor the integer $N$ in deterministic logarithmic time complexity $O(\log(N)^c)$, $c > 0$ constant. If this observation is valid, the difference of squares method is probably the first integer factorization algorithm of deterministic exponential time $O(N^{1/4})$. The fastest, deterministic, and unconditionally proven integer factorization algorithms in the literature have running time complexities $O(N^{1/4})$, see [CP, p. 238]. Many algorithms of deterministic exponential time complexities $O(N^{1/4})$ have been discovered, see [CP], [CE], [MP], et cetera. A relevant integer factorization algorithm in this class is the following.

**Theorem 3.** ([CR])  If the least (or most) significant $(\log N)/4$ bits of a prime factor $p$ or $q$ of the integer $N = pq$, $p < q < 2p$, are given, then it can be decomposed in deterministic logarithmic time complexity $O(\log(N)^c)$, $c > 0$ constant.

This result is equivalent to an integer factorization algorithm of running time complexity $O(N^{1/4}\log(N)^c)$. For example, with respect to this algorithm (Theorem 3), a 1024-bit balanced integer has a running time complexity of approximately 1024/4 = 256 bits. Experimental data for this algorithm are compiled in [CO] and similar references.

**Proposition 4.**  Let $c > 0$ be a constant and let $p$ and $q$ be prime numbers. Then the subset of balanced integers $N = pq$, $p < q < cp$, is of cardinality

$$\#\{\ N = pq \leq x : p < q < cp\ \} = c_1 x / \log^2 x + O(x / \log^3 x)\,.$$

Proof: For $x > x_0$, consider the primes $p \in [c^{-1}N^{1/2}, N^{1/2}]$ and $q \in [N^{1/2}, cN^{1/2}]$. By the standard version of the Prime Number Theorem $\pi(x) = x/\log x + O(xe^{-c_0\sqrt{\log x}})$, (due to delaValle Poussin), it is easy to show that the short interval $[x, x + O(xe^{-c_0\sqrt{\log x}})]$ contains primes. In particular, the number of primes in the previous short intervals are $\pi(x^{1/2}) - \pi(c^{-1}x^{1/2})$ and $\pi(cx^{1/2}) - \pi(x^{1/2})$, respectively. Hence, the subset of balanced integers $N = pq$, $p < q < cp$, has a cardinality of

$$[\pi(x^{1/2}) - \pi(c^{-1}x^{1/2})][\pi(cx^{1/2}) - \pi(x^{1/2})] = c_1 x / \log^2 x + O(x / \log^3 x)\,,$$

where $c_0, c_1 > 0$ are constants.  ∎

For a constant $c > 0$, this is a simple and straightforward counting argument, and it actually follows a calculation similar to the Bertrand postulate. Since the density of primes in the short interval $(x, x + y]$ is known to satisfies the asymptotic expression $\pi(x + y) - \pi(x) = y/\log x + O(y/\log^2 x)$ for $y \geq x^{7/12}$, which is significantly better than the standard version of the Prime Number Theorem, the restriction $c > 0$ constant can be to $c(x)$ = slowly increasing function of $x$ as $x \to \infty$. A different way of proving this result appears in [DM].

In general, the number of squarefree integers with $k$ prime factors is given by

$$\pi_k(x) = \frac{6}{\pi^2(k-1)!} \frac{x(\log\log x)^{k-1}}{\log x} + O(\frac{x(\log\log x)^{k-2}}{\log x})\,,$$

this due to Gauss and landau, see [DN, Vol. I, p. 438].

## 3. DETERMINISTIC EXPONENTIAL TIME $O(N^{1/6})$

For a pair of fixed parameters $0 < \alpha < 1 < \beta$, let $p$ and $q$ be prime numbers such that $\sqrt{\alpha N} < p < \sqrt{N}$, $\sqrt{N} < q < \sqrt{\beta N}$, and let $\gamma\sqrt{N} = (\sqrt{\alpha N} + \sqrt{\beta N})/2$ be the arithmetic mean of the interval. Shifting the symmetric center $\sqrt{N}$ of the Fermat difference of squares method (Theorem 2) or the Coopersmith algorithm (Theorem 3) to the arithmetic mean center $\gamma\sqrt{N}$ of the factors, or to a pair of distinct centers $\sqrt{\alpha N}$ and $\sqrt{\beta N}$, $\alpha, \beta \neq 1$, can augment the ranges of the integers that can be decomposed and reduces the time complexities of both the difference of squares method and the Coopersmith algorithm, respectively.

**Theorem 5.** Let $\alpha > 0$ and $\beta > 0$ be small fixed parameters such that $\alpha\beta - 1 \neq 0$. If the factors of integer $N = pq$ satisfy $|\sqrt{\alpha N} - p| = O(N^{1/3}\log^c N)$ and $|\sqrt{\beta N} - q| = O(N^{1/3}\log^c N)$, then $N$ can be decomposed in deterministic logarithmic time complexity $O(\log(N)^c)$, $c > 0$ constant.

Proof: Examine the product $N = (\sqrt{\alpha N} - x)(\sqrt{\beta N} - y)$. The corresponding polynomial is

$$f(x, y) = xy - \sqrt{\alpha N}\,x - \sqrt{\beta N}\,y + (\alpha\beta - 1)N, \tag{6}$$

where $0 \leq |x|, |y| < 2N^{1/3}$. By construction, this is an irreducible polynomial over the integers. The bounds $X = 2N^{1/3}$ and $Y = 2N^{1/3}$ and the constraint $\alpha\beta - 1 \neq 0$ imply that the height $W = \| f(xX, yY) \|_\infty$ of the polynomial $f(xX, yY)$ is $\| f(xX, yY) \|_\infty = (\alpha\beta - 1)N$. Further, since the height inequality $XY < W^{2/3} = ((\alpha\beta - 1)N)^{2/3} \approx N^{2/3}$ is satisfied, the small integer roots $(x_0, y_0)$ such that $0 \leq |x_0| < X = 2N^{1/3}$ and $0 \leq |y_0| < Y = 2N^{1/3}$ can be determined in deterministic logarithmic time $O(\log(N)^c)$, using lattice reduction methods, see Theorem 11. ∎

Some of the possible values of the parameters $\alpha$ and $\beta$, of interest, are approximately $\alpha \approx 1/2$, $\beta \approx 2$.

**Example 2.** For $\alpha = 1/1.85$, $\beta = 2.15$, the polynomial is $f(x, y) = xy - \sqrt{.541N}\,x - \sqrt{2.150N}\,y + 0.162N$.

**Example 3.** For $\alpha = 1/2$, $\beta = 2$, and the arithmetic mean $(\sqrt{\alpha N} + \sqrt{\beta N})/2 = \gamma\sqrt{N}$ of the prime factors, the corresponding polynomial is $f(x, y) = (\gamma\sqrt{N} - x)(\gamma\sqrt{N} + y) - N = xy + 1.061\sqrt{N}\,x - 1.061\sqrt{N}\,y + 0.125N$. The factors $p = \gamma\sqrt{N} - x$ and $q = \gamma\sqrt{N} + y$ of the integers $N$ such that $0 \leq |x| < N^{1/3}$ and $0 \leq |y| < N^{1/3}$ can be determined using the lattice reduction method. In addition, the difference of squares method started at $x_k = \gamma\sqrt{N} + k$, $k \geq 0$, can be used to factor $N$ provided that the starting point $2\gamma\sqrt{N} \leq p + q$.

At the values $\alpha = \beta = 1$, the arithmetic mean of the prime factors is $(\sqrt{N} + \sqrt{N})/2 = \sqrt{N}$, so this algorithm (Theorem 5) reduces to the factoring methods of Theorems 2 and 3.

Extending the range from $0 \leq |x|, |y| < 2N^{1/3}$ to $0 \leq |x|, |y| < 2N^{1/2}$ transforms the algorithm (Theorem 5) into a general purpose integer factorization algorithm of deterministic exponential time complexity $O(N^{1/6}\log^c N)$.

**Corollary 6.** Given the least (or most) significant $(\log N)/6$ bits of a prime factor $p$ or $q$ of a large integer $N = pq$, $p < q < 2p$, the integer $N$ can be decomposed in deterministic logarithmic time $O(\log^c N)$, $c > 0$ constant.

Proof: Assume that $x_0, y_0$ are the given least significant $(\log N)/6$ bits of the prime factors, ($y_0$ is computed via the congruence $N_0 \equiv -x_0 y_0 \bmod N^{1/6}$, where $N_0 = [(\gamma^2 - 1)N]$, $\gamma^2 - 1 \neq 0$, and the given $x_0$ in $p = \gamma N^{1/2} + xN^{1/6} + x_0$, and put $N = (\gamma N^{1/2} + xN^{1/6} + x_0)(\gamma N^{1/2} + yN^{1/6} + y_0)$, where $0 \leq |x|, |y| < N^{1/3}$. By construction, the corresponding polynomial

$$f(x, y) = xyN^{1/3} + (\gamma N^{1/2} + y_0)N^{1/6}x + (\gamma N^{1/2} + x_0)N^{1/6}y + (x_0 + y_0)\gamma N^{1/2} + (\gamma^2 - 1)N + x_0 y_0, \tag{7}$$

where $\gamma = (\sqrt{2} + \sqrt{1/2})/2 = 3\sqrt{2}/4$ is the arithmetic mean, is an irreducible polynomial over the integers. For the bounds $X = N^{1/3}$ and $Y = N^{1/3}$, the height $W = \| f(xX, yY) \|_\infty$ of the polynomial $f(xX, yY)$ is $\| f(xX, yY) \|_\infty = N$. Thus, using lattice reduction methods, the small integer roots $0 \leq |x| < X = N^{1/3}$ and $0 \leq |y| < Y = N^{1/3}$ can be determined in deterministic logarithmic time, see Theorem 11. ∎

As an application, consider the time complexities of factoring 686-bit and 1024-bit integers $N = pq$, $p < q < 2p$. Relative to this algorithm (Corollary 6), these integers have 686/12 ≈ 60 bits and 1024/12 ≈ 86 bits time complexities, respectively. In other words, to factor 686-bit and 1024-bit balanced integers, it is sufficient to conduct brute force searches for 60 consecutive bits and 86 consecutive bits of $p$ or $q$ (or given 60 consecutive bits and 86 consecutive bits of $p$ or $q$, respectively), these integers can be factored using the lattice reduction methods unconditionally. The current estimates of the costs of factoring these integers, using the number fields sieve algorithm are discussed in [BK] and [KZ].

## 4. ELEMENTARY CONCEPTS IN SYSTEMS OF POLYNOMIALS EQUATIONS

The research in the theory of systems of polynomials equations covers a wide spectrum of the mathematical literature, from algebraic geometry to numerical analysis and beyond. The simplest case of polynomial equations of a single variable has a vast and active literature. The algorithms for determining real/complex roots of a polynomial equation $f(x) = 0$ over the integers are efficient, but the algorithms for determining the solutions of the modular case $f(x) \equiv 0 \bmod N$, where $N$ is an arbitrary integer, are not completely effective yet. The theory of quadratic polynomial equations $f(x,y) = 0$ over the integers is not completely understood. For example, determining whether or not $ax^2 + by + c = 0$ has an integer solution is classified as an NP-complete problem. However, there are various partial results obtained by lattice reduction methods and other techniques. The applications of lattice reduction methods to the theory of polynomial equations and its applications to cryptography are considered in fine details in [VE], [CR], [HR], [BM], [CO], [JZ], [LH] and others.

Several relevant results from the theory of polynomials equations are included in this section to set the notations and to provide some background materials. A comprehensive introduction to lattice reduction methods and its applications to polynomial equations is given in [JZ, Chapter 3] and similar sources. The evolving analysis on a few specific polynomial equations of three variables is given in [BA].

***Definition 7.*** A subset of polynomials $f_1(x_1, \ldots, x_n), \ldots, f_m(x_1, \ldots, x_n) \in \mathbb{Z}[x_1, \ldots, x_n]$ are called algebraically independent if and only if the relation $P(f_1, \ldots, f_m) = 0$ implies that $P(t_1, \ldots, t_m) \in \mathbb{Z}[t_1, \ldots, t_m]$ is the zero polynomial.

The concept of algebraic independence is covered in finer details in [CX], [ML] and related literature. A crucial property of pairwise algebraically independent polynomials is the vanishing/nonvanishing behavior of the resultant. For pairwise algebraically independent polynomials, the (first level) resultants $R(f, g, x_i) \neq 0$ with respect to any of the variables $x_i$, $0 \leq i \leq n$, are nonvanishing functions of $n - 1$ variables $x_1, \ldots, x_n$. But the (second level) resultants $R(R(f, g, x_i), R(f, g, x_i), x_j) = 0$ with respect to any pair of variables $x_i \neq x_j$ vanish.

The euclidean norm $\| f \|_2$ and the supremum norm $\| f \|_\infty$ of a polynomial $f(x, y, z) = \sum_{0 \leq i,j,k \leq d} a_{i,j,k} x^i y^j z^k \in \mathbb{Z}[x, y, z]$ over the integers are defined by

$$\| f \|_2^2 = \sum_{0 \leq i,j,k \leq d} |a_{i,j,k}|^2 \quad \text{and} \quad \| f \|_\infty = \max\{ |a_{i,j,k}| \}, \tag{8}$$

respectively.

**Theorem 8.** ([ST]) Let $f(x_1, \ldots, x_n)$, $g(x_1, \ldots, x_n) \in \mathbb{Z}[x_1, \ldots, x_n]$ be nonzero polynomials of maximum degrees $d > 0$ in each variable separately, such that $g(x_1, \ldots, x_n)$ is a multiple of $f(x_1, \ldots, x_n)$. Then

$$\| g \|_2 \geq 2^{-(d+1)^n + 1} \| f \|_\infty .$$

This is a straightforward derivation from the earlier result in [MT] for polynomials of a single variable.

Let $\{ U_1, U_2, \ldots, U_d \} \subset \mathbb{R}^n$ be a subset of vectors in the $n$-dimensional real vector space $\mathbb{R}^n$. A discrete lattice $L = \{ a_1 U_1 + a_2 U_2 + \cdots + a_d U_d : a_i \in \mathbb{R} \} \subset \mathbb{R}^n$ is a discrete subset of vector space $\mathbb{R}^n$. Define the sequence of numbers $\lambda_1(L) \leq \lambda_2(L) \leq \cdots \leq \lambda_n(L)$ as the minimum norms $\| V_i \|_2$ of any subset of vectors in a lattice $L$. A result from the geometry of numbers specifies the possible upper bound of the norms of any subset of $d$ vectors, $d \leq n$.

**Theorem 9.** (Minkowski) Let $L$ be a lattice of dimension $n \geq 1$. Then:
(i) The norm of some vector $V \in L$ satisfies the inequality $\| V \|_2 \leq \sqrt{n} \det(L)^{1/n}$.
(ii) The product of the norms of the $d$ smallest vectors $V_i$ satisfies the inequality $\prod_{i \leq d} \| V_i \|_2 \leq \prod_{i \leq d} \lambda_i(L) \leq \gamma_n^{d/2} \det(L)^{d/n}$, where $n/(2e\pi) \leq \gamma_n \leq n/(e\pi + o(1))$ is the $n$th Hermite constant, and $e = 2.71\ldots$ is the base of the logarithmus naturalis.

**Theorem 10.** ([LH]) Let $L$ be a lattice of dimension $n \geq 1$. Then there is an algorithm that generates a reduced basis $\{ V_1, V_2, \ldots, V_n \}$ with the property

$$\| V_1 \|_2 \leq \| V_2 \|_2 \leq \cdots \leq \| V_n \|_2 \leq 2^{n(n-1)/(4(n+1-i))} \det(L)^{1/(n+1-i)}$$

for $0 \leq i \leq n$, in deterministic logarithmic time $O(\log^c N)$, $c > 0$ constant.

On average, the norm of a small vector in a lattice is approximated by the geometric mean, that is, $\| V_i \|_2 \leq \gamma_n^{1/2} \det(L)^{1/n}$. By comparison, the LLL algorithm determines a short vector in the lattice within a factor of $2^{n-1}$ of the geometric mean, that is, $\| V_i \|_2 \leq 2^{n-1} \det(L)^{1/n}$, and this is accomplished in effective time complexity.

Employing lattice reduction methods, several results for the polynomials $f(x, y) = \sum_{0 \leq i, j \leq d} a_{i,j} x^i y^j$ and $f(x, y, z) = \sum_{0 \leq i, j, l \leq d} a_{i,j,l} x^i y^j z^l$ of two and three variables have been unconditionally proven.

**Theorem 11.** ([CR]) Let $f(x, y) \in \mathbb{Z}[x, y]$ be an irreducible polynomial of maximum degree $d$, and let $(x_0, y_0)$ be a root of $f(x, y) = 0$, such that $0 \leq |x_0| \leq X$, $0 \leq |y_0| \leq Y$. The height of the polynomial $f(xX, yY)$ is defined by $W = \| f(xX, yY) \|_\infty = \max \{ | a_{i,j} X^i Y^j | : 0 \leq i, j \leq d \}$.
i) If $XY < W^{2/(3d)}$, then the roots $(x_0, y_0)$ can be determined in deterministic logarithmic time $O(\log^c N)$, $c > 0$ constant.
ii) If $XY < W^{1/d}$, and the total degree of the polynomial satisfies $0 \leq i + j \leq d$, then the roots $(x_0, y_0)$ can be determined in deterministic logarithmic time $O(\log^c N)$, $c > 0$ constant.

The detailed heuristic analysis of the lattices for $f(x, y, z) = c_0 xy + c_1 x + c_2 y + c_3 z + c_4 \in \mathbb{Z}[x, y, z]$, and a few other polynomials of three variables and related applications, appear in [ER, p. 8], and [JZ, p. 66]. Practical applications also appear in [JM] and [HK, p. 11].

**Theorem 12.** ([ER], [JZ]) Let $f(x, y, z) = c_0xy + c_1x + c_2y + c_3z + c_4$ be an irreducible polynomial, and let $(x_0, y_0, z_0)$ be a root of $f(x, y, z) = 0$ such that $0 \leq |x_0| \leq X$, $0 \leq |y_0| \leq Y$, $0 \leq |z_0| \leq Z$. Suppose that the lattice reduction algorithm can effectively generate an algebraically independent subset of polynomials { $f(x, y, z)$, $f_1(x, y, z), f_2(x, y, z)$ } contingent on the constraint

$$X^{3+3\tau}Y^{3+6\tau+3\tau^2}Z^{2+3\tau} < W^{2+3\tau-\varepsilon}, \tag{9}$$

where $W = \| f(xX, yY, zZ) \|_\infty$ is the height, $\varepsilon > 0$ is a small number, and $\tau > 0$ is a lattice parameter. Then the roots of the system of equations $f(x, y, z) = 0$, $f_1(x, y, z) = 0$, $f_2(x, y, z) = 0$ can be determined in deterministic logarithmic time complexity $O(\log^c N)$, $c > 0$ constant.

## 5. HEURISTIC LOGARITHMIC TIME $O(\log(N)^c)$

Although the lattice reduction methods analysis for polynomials of three or more variables is heuristic, and proved in some special cases, see [BA], the numerical experiments are very effective, see [CR], [CO], [BM], [ER], [JM], [JZ] and similar literature. In light of these numerical results, it is expected that the heuristic integer factorization algorithm developed here is effective.

**Theorem 13.** Assume that the lattice reduction method can generate two algebraically independent polynomials, then the prime factors of the integer $N = pq$, $p < q < 2p$, can be determined in deterministic logarithmic time $O(\log^c N)$, $c > 0$ constant.

Proof: Let $p_0, q_0 = O(N^{1/2-\alpha})$, $r_0 \neq 0$, and $m_0 = N^{1+\beta}$ be fixed integer parameters such that $\alpha > 0$, $\beta > 0$ are small numbers and $\gcd(m_0, r_0) = 1$. Now consider the identities

$$r_0 + m_0p = r_0 + m_0(p_0 + x) \quad \text{and} \quad r_0 + m_0q = r_0 + m_0(q_0 + y), \tag{10}$$

where $0 \leq |x|, |y| < N^{1/2}$. Taking the product of these terms leads to

$$r_0^2 + r_0m_0(p+q) + m_0^2N = r_0^2 + r_0m_0(x+p_0) + r_0m_0(y+q_0) + m_0^2(x+p_0)(y+q_0). \tag{11}$$

Replace $z = p + q$, and rearrange it to obtain the polynomial

$$f_0(x, y, z) = m_0xy + (r_0 + m_0q_0)x + (r_0 + m_0p_0)y - r_0z + m_0p_0q_0 + r_0(p_0 + q_0) - m_0N. \tag{12}$$

The integer parameters are properly selected to construct an irreducible polynomial $f_0(x, y, z)$ over the integers of sufficiently large height. For the bounds $X = N^{1/2}$, $Y = N^{1/2}$, and $Z = 3N^{1/2}$, the height $W = \| f_0(xX, yY, zZ) \|_\infty$ of the polynomial $f_0(xX, yY, zZ)$ is given by

$$\| f_0(xX, yY, zZ) \|_\infty \geq | m_0p_0q_0 + r_0(p_0 + q_0) - m_0N | \geq m_0N = N^{2+\beta}. \tag{13}$$

Substituting the limits $X = N^{1/2}$, $Y = N^{1/2}$, and $Z = 3N^{1/2}$ in the *determinant/height inequality* (9), returns

$$cN^{(3+3\tau)/2}N^{(3+6\tau+3\tau^2)/2}N^{(2+3\tau)/2} = cN^{4+6\tau+3\tau^2/2} < N^{(2+\beta)(2+3\tau-\varepsilon)}, \tag{14}$$

where $c > 0$ is a small constant, $\varepsilon > 0$, and $\tau \geq 0$ are lattice parameters, see [JZ, p. 66]. This requires

$(3\tau^2/2 + 2\varepsilon)/(2 + 3\tau - \varepsilon) < \beta$. But, since there is no restriction on the small parameter $\beta > 0$, this inequality holds for suitable values of the parameters $\beta > 0$, $\varepsilon > 0$, $\tau \geq 0$. Under these conditions, and by hypothesis, the lattice reduction algorithm can produce a triple of algebraically independent polynomials

$$\{ \; f_0(x, y, z), \; f_1(x, y, z), \; f_2(x, y, z) \; \}, \tag{15}$$

see Theorem 12. Hence, the system of equations

$$f_0(x, y, z) = 0, \; f_1(x, y, z) = 0, \; f_2(x, y, z) = 0 \tag{16}$$

can be solved for the small integer solutions $0 \leq |x_0| < N^{1/2}$, $|y_0| < N^{1/2}$, $|z_0| < 3N^{1/2}$ in deterministic logarithmic time $O(\log^c N)$, $c > 0$ constant. Last, but not least, the prime factors are recovered using the expressions $p = p_0 + x_0$ and $q = q_0 + y_0$. ∎

The integer parameters $p_0$, $q_0 = O(N^{1/2-\alpha})$, $r_0 \neq 0$, and $m_0 = N^{1+\beta}$ facilitates the construction of irreducible polynomials of large height and various shapes suitable for the lattice reduction algorithms. These parameters are easy to adjust to optimize the polynomials, the lattices, and the algorithms. Indeed, almost any random pair $p_0$, $q_0 = O(N^{1/2-\alpha})$, and a pair of fixed integers $m_0 = N^{1+\beta}$, and $r_0 \neq 0$ such that $\gcd(m_0, r_0) = 1$, with small numbers $\alpha > 0$, $\beta > 0$, can achieve these requirements. A possible implementation of this algorithm is sketched below.

**Algorithm I.**
Input: $N = pq$, $p < q < 2p$.
Output: $p$, $q$.

1. *Irreducible Polynomial Routine.*
1.1 Put $p_0 = 2[.492343N^{.378549}] + 1$, $q_0 = 2[.649287N^{.487532}] + 1$, $m_0 = 2[841.013799N^2] + 1$, and fix an integer $r_0 = \pm 2^k$, $k \geq 0$, such that $\gcd(m_0, r_0) = 1$, the bracket $[\;x\;]$ denotes the greatest integer function. This is one of the possible selection criteria, these choices ensure that $\gcd(c_0, c_1, c_2, c_3, c_4) = 1$, the polynomial is irreducible, and that it has sufficiently large height $\| f_0(xX, yY, zZ) \|_\infty \geq 1682N^3$.
1.2 Put $f_0(x, y, z) = c_0xy + c_1x + c_2y + c_3z + c_4$
$$= m_0xy + (r_0 + m_0q_0)x + (r_0 + m_0p_0)y - r_0z + m_0p_0q_0 + r_0(p_0 + q_0) - m_0N.$$

2. *Basis Routine.*
2.1 Call a lattice reduction routine to compute the subset $\{ \; f_0(x, y, z), \; f_1(x, y, z), \; f_2(x, y, z) \; \}$.

3. *Roots Routine.*
3.1 Call a resultant and polynomial root routines (or Gröbner basis routine) to compute the small roots $(x_0, y_0, z_0)$ such that $0 \leq |x_0| < N^{1/2}$, $|y_0| < N^{1/2}$, $|z_0| < 3N^{1/2}$ of the system of polynomial equations $f_0(x, y, z) = 0$, $f_1(x, y, z) = 0$, $f_2(x, y, z) = 0$.

4. *Factors Routine.*
4.1 Return $p = p_0 + x_0$ and $q = q_0 + y_0$.

## 6. UNCONDITIONAL RESULT

A bench mark for the time complexity of the integer factorization problem was established in [CL], but there was no claim on the practicality of the algorithm. The rudimentary details of a different, and more practical and unconditional proof are sketched on this Section. There are a few approaches to an unconditional proof of the heuristic algorithm of the previous section.

(1) Lattice reduction methods.          (2) Purely algebraic methods.
(3) Purely number theoretical methods.

A combination of lattice reduction methods and algebraic methods will be utilized here. These tools seem to be sufficient and practical for numerical experiments. A purely algebraic proof in terms of ideals and varieties is also very interesting but perhaps a bit more difficult, see [BA], [CX], [ML] etc. Similarly, a purely number theoretical proof seems to be difficult too.

It is quite easy to construct a subset of two or more pairwise algebraically independent polynomials $\{v_0(x, y, z), v_1(x, y, z)\}$ by means of elementary polynomial operations, see Definition 7. But the third polynomial needed to complete the subset of at least three triplewise algebraically independent polynomials $\{v_0(x, y, z), v_1(x, y, z), v_2(x, y, z)\}$ seems to be constructible only by nonelementary polynomial operations, such as lattice reduction techniques, and other means.

The basic generating relations $N = pq$, and $z = p + q$ for the integer factorization problem and elementary polynomial operations are utilized to construct the subset of polynomials

$$u_0(X, Y, z) = X + Y - z, \tag{17}$$

$$u_1(X, Y, z) = X^2 - zX + N, \qquad\qquad u_2(X, Y, z) = Y^2 - zY + N,$$

$$u_3(X, Y, z) = X^4 - (z^2 - 2N)X^2 + N^2, \qquad u_4(X, Y, z) = Y^4 - (z^2 - 2N)Y^2 + N^2,$$

$$u_5(X, Y, z) = X^6 - (z^3 - 3Nz)X^3 + N^3, \qquad u_6(X, Y, z) = Y^6 - (z^3 - 3Nz)Y^3 + N^3,$$

$$u_7(X, Y, z) = X^8 - (z^4 - 4Nz^2 + 2N^2)X^4 + N^4, \qquad u_8(X, Y, z) = Y^8 - (z^4 - 4Nz^2 + 2N^2)Y^4 + N^4,$$

$$\cdots \qquad\qquad\qquad\qquad\qquad\qquad \cdots$$

$$f_0(x, y, z) = m_0 xy + (r_0 + m_0 q_0)x + (r_0 + m_0 p_0)y - r_0 z + m_0 p_0 q_0 + r_0(p_0 + q_0) - m_0 N,$$

$$f_3(x, y, z) = m_1 xy + (r_1 + m_1 q_0)x + (r_1 + m_1 p_0)y - r_1 z + m_1 p_0 q_0 + r_1(p_0 + q_0) - m_1 N,$$

et cetera. The variables are linked by the linear change of variables $X = x + p_0 = p$, $Y = y + q_0 = q$. The parameters of the polynomials $f_0(x, y, z)$ and $f_3(x, y, z)$ above satisfy $\gcd(m_0, m_1, r_0, r_1) = 1$, see (10), (11) and (12) for more details on the derivation of $f_0(x, y, z)$ and $f_3(x, y, z)$.

Since the (first level) resultants

$$R(u_i(X, Y, z), u_j(X, Y, z), X) \neq 0, \ R(u_i(X, Y, z), u_j(X, Y, z), Y) \neq 0, \ R(u_i(X, Y, z), u_j(X, Y, z), z) \neq 0 \tag{18}$$

do not vanish for most $i \neq j$, this list contains various pairs of pairwise algebraically independent polynomials. In addition, the elementary resultant operation can be used to reproduce new polynomials, for example,

$$R(u_0(X, Y, z), u_3(X, Y, z), z) = -2X^3 Y + (2N - Y^2)X^2 + N^2 \neq 0, \tag{19}$$

$$R(u_0(X, Y, z), u_5(X, Y, z), z) = -3X^5 Y + (3N - 3Y^2)X^4 + (3NY - Y^3)X^3 + N^3 \neq 0,$$

and so on. But the next (second level) of resultants

$$R(R(u_i(X, Y, z), u_j(X, Y, z), z), R(u_k(X, Y, z), u_l(X, Y, z), z), Y) = 0 \tag{20}$$

for $0 \leq i, j, k, l$, vanish because these polynomials are not triplewise algebraically independent. More precisely, the definition of algebraic independence, $P(g_1, g_2, g_3) = 0$ for some nonzero polynomial $P(t_1, t_2, t_3) \in \mathbb{Z}[t_1, t_2, t_3]$ and any three polynomials $g_1$, $g_2$, $g_3$ from the list of polynomials (17).

The aim of this section is to utilize nonelementary polynomial operations to complete the subset (17) into a subset of triplewise algebraically independent polynomials. This task is the core of the proof of the next result.

***Theorem* 14.** The integer factorization problem has deterministic logarithmic time complexity $O(\log^c N)$, $c > 0$ constant.

Proof: Consider the polynomials

$$f_0(x, y, z) = m_0 xy + (r_0 + m_0 q_0)x + (r_0 + m_0 p_0)y - r_0 z + m_0 p_0 q_0 + r_0(p_0 + q_0) - m_0 N , \tag{21}$$
$$f_3(x, y, z) = m_1 xy + (r_1 + m_1 q_0)x + (r_1 + m_1 p_0)y - r_1 z + m_1 p_0 q_0 + r_1(p_0 + q_0) - m_1 N ,$$

where $\gcd(m_0, m_1, r_0, r_1) = 1$, and $m_0 \geq N^{1+\beta}$, $m_1 \geq N^{1+\beta}$, see (10), (11), and (12) for details on the construction. These irreducible polynomials $f_0(x, y, z)$ and $f_3(x, y, z)$ share a common root $(x_0, y_0, z_0)$, such that $p = p_0 + x_0$, $q = q_0 + y_0$ and $p + q = z_0 + p_0 + q_0$. Furthermore, by hand or machine calculations, it can be shown that the subset of polynomials $\{ f_0, f_3 \}$ are pairwise algebraically independent.

The bounds are $X = N^{1/2}$, $Y = N^{1/2}$, and $Z = 3N^{1/2}$, so that the small roots $(x, y, z)$ of interest are in the ranges $0 \leq |x| \leq X = N^{1/2}$, $0 \leq |y| \leq Y = N^{1/2}$, and $0 \leq |z| \leq Z = 3N^{1/2}$.

Now utilize the irreducible $f_0(x, y, z)$ (or $f_3(x, y, z)$) to construct a lattice basis and employ a lattice reduction algorithm to generate the pair of polynomials

$$f_1(x, y, z) = \sum_{0 \leq i,j,k \leq d} a_{i,j,k} x^i y^j z^k , \tag{22}$$
$$f_2(x, y, z) = \sum_{0 \leq i,j,k \leq d} b_{i,j,k} x^i y^j z^k .$$

By Theorem 12, the subset of polynomials $\{ f_0, f_1, f_2 \}$ are pairwise algebraically independent but unknown a priori whether or not these are triplewise algebraically independent, see Definition 7.

Next, to show that the subset $\{ f_0, f_1, f_2, f_3 \}$ contains three algebraically independent polynomials, assume that either $f_1(x, y, z)$ or $f_2(x, y, z)$ is a multiple of $f_3(x, y, z)$. The euclidean norms of the polynomials $f_1(x, y, z)$ and $f_2(x, y, z)$ satisfy the inequalities

$$\| f_1 \|_2 \leq \| f_2 \|_2 \leq 2^{n(n-1)/(4(n+1-i))} \det(L)^{1/(n+1-i)} , \tag{23}$$

where $n \geq 1$ is the dimension of the lattice, and $0 \leq i \leq n$. And the supremum norms of the polynomials $f_0(x, y, z)$ and $f_3(x, y, z)$ satisfy the inequalities

$$2^{n(n-1)/(4(n+1-i))} \det(L)^{1/(n+1-i)} \leq \| f_0 \|_\infty \leq \| f_3 \|_\infty . \tag{24}$$

In addition, the reduced lattice basis $\{ V_1, V_2, \ldots, V_n \}$ satisfies the inequalities

$$\| V_1 \|_2 \leq \| V_2 \|_2 \leq \cdots \leq \| V_n \|_2 \leq 2^{n(n-1)/(4(n+1-i))} \det(L)^{1/(n+1-i)} \leq \| f_0 \|_\infty \tag{25}$$

for some $n + 1 - i = n - 1 \geq 2$, see Theorem 10. These, in turn, imply that the euclidean norms of the polynomials $f_1(x, y, z)$ and $f_2(x, y, z)$ satisfy

$$\|\,f_1\,\|_2 \le \|\,f_2\,\|_2 \le 2^{n(n-1)/(4(n-1))}\det(L)^{1/(n-1)} \le \|\,f_3\,\|_\infty. \tag{26}$$

Therefore, by Theorem 8, neither the polynomial $f_1(x,\,y,\,z)$ nor $f_2(x,\,y,\,z)$ is a multiple of the irreducible polynomial $f_3(x,\,y,\,z)$, confer [CO] for similar argument in the case of two variables. Consequently, the subsets of polynomials

$$\{\,f_0,\,f_1\,\},\ \{\,f_0,\,f_2\,\},\ \{\,f_0,\,f_3\,\},\ \{\,f_1,\,f_3\,\},\ \{\,f_2,\,f_3\,\} \tag{27}$$

are pairwise algebraically independent. The algebraic independence status of $\{\,f_1,\ f_2\,\}$ is unknown a priori. Therefore, there exists a subset of three algebraically independent polynomials, e. g.,

$$\{\,f_0(x,y,z),\ f_1(x,y,z),\ f_i(x,y,z)\ \}, \tag{28}$$

where $i \in \{\,2,\,3\,\}$. Specifically, the Gröbner basis $\{\,g_0(x,y,z), g_1(x,y,z),...,\,g_m(x,y,z)\ \}$ of the set

$$\{\,f_0(x,y,z),\,f_1(x,y,z),\,f_2(x,y,z),\,f_3(x,y,z)\ \} \tag{29}$$

contains a polynomial of a single variable $g_i(x,\,y,\,z) = v(x)$ or $v(y)$ or $v(z)$. Equivalently, the system of equations

$$f_0(x,y,z) = 0,\,f_1(x,y,z) = 0,\,f_2(x,y,z) = 0,\,f_3(x,y,z) = 0 \tag{30}$$

can be solved for the small roots $(x_0, y_0, z_0)$. 

Quod Erat Demonstrandum.

Note that the subset of polynomials

$$\{\,f_0(x,y,z),\,f_1(x,y,z),\,f_2(x,y,z),f_4(x,y,z) = x + y - z + p_0 + q_0\ \} \tag{31}$$

and other simple combinations of polynomials from the list (17) have simpler resultants and Gröbner bases calculations. But, the proof of Theorem 14 appears to be more difficult. In particular, the step about proving that neither the polynomial $f_1(x,\,y,\,z)$ nor $f_2(x,\,y,\,z)$ is a multiple of $f_4(x,\,y,\,z) = x + y - z + p_0 + q_0$ is not as simple as above. However, in practice, the simpler subset

$$\{\,f_0(x,y,z),\,f_1(x,y,z),\,f_4(x,y,z)\ \} \tag{32}$$

and its corresponding system of equations

$$f_0(x,y,z) = 0,\,f_1(x,y,z) = 0,\,f_4(x,y,z) = 0 \tag{33}$$

should be sufficient to handle the integer factorization problem in effective running time complexity.